%%%%%%%%%%%%%%%%%%%%%%%%%%%%%%%%%%%%%%%%%%%%%%%%%%%%%%%%%%%%%%%%%%%%%%%%%%%
%% Boyer, S.; Zhang, X.
%% 
%% Exceptional surgery on knots
%% 
%% Let $M$ be an irreducible, compact, connected, orientable 3-manifold 
%%   whose boundary is a torus. We show that if $M$ is hyperbolic, then it 
%%   admits at most six finite/cyclic fillings of maximal distance 5. 
%%   Further, the distance of a finite/cyclic filling to a cyclic filling 
%%   is at most 2. If $M$ has a non-boundary-parallel, incompressible torus
%%   and is not a generalized 1-iterated torus knot complement, then there 
%%   are at most three finite/cyclic fillings of maximal distance 1. 
%%   Further, if $M$ has a non-boundary-parallel, incompressible torus and 
%%   is not a generalized 1- or 2-iterated torus knot complement and if $M$
%%   admits a cyclic filling of odd order, then $M$ does not admit any other
%%   finite/cyclic filling. Relations between finite/cyclic fillings and 
%%   other exceptional fillings are also discussed.
%% 
%% publ:  Bull. Amer. Math. Soc. (N.S.) 31(1994) no. 2
%% pp:    197-203
%% type:  Research Announcement        markup: amstex    file size: 26K
%% contact: boyer@math.uqam.ca
%% 
%% copyright: American Math. Society copyright; see end of article
%% 
%% Include files necessary for this article: bull-ppt.tex
%% 
%%%%%%%%%%%%%%%%%%%%%%%%%%%%%%%%%%%%%%%%%%%%%%%%%%%%%%%%%%%%%%%%%%%%%%%%%%%
\input amstex 
\documentstyle{amsppt}
\input bull-ppt
\keyedby{bull516/kmt}

\topmatter
\cvol{31}
\cvolyear{1994}
\cmonth{October}
\cyear{1994}
\cvolno{2}
\cpgs{197-203}
%\ratitle
\title Exceptional surgery on knots \endtitle
\author S. Boyer and X. Zhang\endauthor
%\shortauthor{}
%\shorttitle{}
\address  D\'epartement de Math\'ematiques et 
d'Informatique, 
Universit\'e du Qu\'ebec \`a Montr\'eal, P.O. 
Box
8888, Station A, Montr\'eal H3C 3P8, Canada\endaddress
\ml\nofrills{\it E-mail address}\,:\enspace
\tt boyer\@math.uqam.ca
\mlpar
{\it E-mail address}\,:\enspace zhang\@math.uqam.ca\endml
%\cu \endcu
\date September 14, 1993\enddate
%\keywords{}
\thanks The first author was partially supported by grants 
NSERC OGP 0009446
and FCAR EQ 3518. The second author was supported by a 
postdoctoral fellowship
from the Centre de Recherches Mathematiques\endthanks
\thanks A version of part of the material of this paper 
was presented by the
first author at the joint AMS-SMM meeting in Merida, 
Mexico, in December 1993\endthanks
\abstract Let $M$ be an irreducible, compact, connected, 
orientable
3-manifold
whose boundary is a torus. We show that if $M$ is 
hyperbolic, then it admits
at most six finite/cyclic fillings of maximal distance 5. 
Further, the
distance of a finite/cyclic filling to a cyclic filling is 
at most 2. If $M$
has a non-boundary-parallel, incompressible torus and is 
not a generalized
1-iterated torus knot complement, then there are at most 
three finite/cyclic
fillings of maximal distance 1. Further, if $M$ has a 
non-boundary-parallel,
incompressible torus and is not a generalized 1- or 
2-iterated torus knot
complement and if $M$ admits a cyclic filling of odd 
order, then $M$ does not
admit any other finite/cyclic filling. Relations between 
finite/cyclic fillings
and other exceptional fillings are also 
discussed.\endabstract
\endtopmatter

\document
\par Let $M$ be a compact, connected, orientable 
3-manifold whose boundary is a
torus. $M$ is called {\it irreducible} if any embedded 
2-sphere in $M$ bounds
an embedded 3-ball. Otherwise, $M$ is called {\it 
reducible}. A {\it slope} on
$\partial M$ is a $\partial M$-isotopy class $r$ of 
essential, simple closed curves. A
slope $r$ is called a {\it boundary} slope if there is a 
properly embedded,
orientable, incompressible, $\partial$-incompressible 
surface $F$ in $M$ such
that $\partial F$ is a nonempty set of parallel simple 
closed curves on
$\partial M$ of slope $r$. The {\it distance} between two 
slopes $r_1$ and
$r_2$, denoted by $\Delta(r_1,r_2)$, is the minimal 
geometric intersection
number between $r_1$ and $r_2$. As $\partial M$ is a 
torus, $\Delta(r_1,r_2)$
may be calculated as the absolute value of the algebraic 
intersection number of
the homology classes carried by $r_1$ and $r_2$.
\par Given a slope $r$ on $\partial M$, a well-defined 
closed 3-manifold $M(r)$ can be
constructed by attaching a solid torus to $M$ via a 
homeomorphism which
identifies a meridian curve of the solid torus to a 
representative curve for
$r$. $M(r)$ is called the $r$-{\it filling} of $M$. Given 
a knot $K$ in a
closed, connected, orientable 3-manifold $W$, with tubular 
neighbourhood $N(K)$,
exterior $M=W-\roman{int}(N(K))$, and slope $r$ on 
$\partial M$, $M(r)$ is also referred
to as the $r$-{\it surgery} of $W$ along $K$.
\par A fundamental result of Wallace [Wa] and Lickorish 
[Li] states that each
closed, orientable 3-manifold results from surgery on some 
link in the
3-sphere. Thus a natural approach to 3-manifold topology 
is to analyze to what
extent various aspects of the topology of a manifold $M$, 
as above, are
inherited by the manifolds $M(r)$. For instance, one could 
try to determine
when a closed essential surface in $M$ becomes inessential 
in some $M(r)$ or
when an irreducible $M$ could produce a reducible $M(r)$. 
An excellent survey
of this topic may be found in [Go1]. Another example of 
some importance arises
when $M$ is a hyperbolic manifold, i.e., $\roman{int}(M)$ 
admits a complete
hyperbolic metric of finite volume. Thurston [Th] has 
shown that in this
situation, all but finitely many of the manifolds $M(r)$ 
are hyperbolic. The
nonhyperbolic slopes include those whose fillings are:
\roster
\item"(a)" manifolds with finite or cyclic fundamental 
groups,
\item"(b)" manifolds which are reducible,
\item"(c)" manifolds which are Seifert-fibred spaces,
\item"(d)" manifolds which admit an incompressible torus.
\endroster
Thurston's geometrication conjecture [Th] predicts that 
the remaining slopes
yield fillings which are hyperbolic.
\par A basic problem then is to describe the set of {\it 
exceptional} slopes on
a torally bounded manifold $M$, i.e., those slopes 
producing nongeneric
fillings. An appropriate description should include an 
upper bound on the
number of such slopes as well as a qualitative measure of 
their relative
positions determined by bounds on their mutual distances. 
We call a slope $r$
on $\partial M$ {\it cyclic} if $M(r)$ has cyclic 
fundamental group; and similarly we
shall refer to finite slopes, reducible slopes, Seifert 
slopes, and essential
torus slopes. It is a remarkable result, the {\it cyclic 
surgery theorem}
[CGLS], that if $M$ is not a Seifert-fibred space, then 
all cyclic slopes on
$\partial M$ have mutual distance no larger than 1, and 
consequently there are at
most three such fillings. Gordon and Luecke obtained 
similar estimates for the
set of reducible slopes [GLu1] and also examined essential 
torus slopes [Go2,
GLu2]. In another direction, Bleiler and Hodgson [BH] 
refined results of Gromov
and Thurston to obtain restrictions on the slopes on a 
hyperbolic $M$ which do
not produce manifolds admitting metrics of strictly 
negative sectional
curvature. We consider the problem of determining the set 
of finite slopes on
$M$ here. Henceforth, we shall use finite/cyclic to mean 
either infinite cyclic
or finite. Standard arguments show that we may take $M$ to 
be irreducible, and
we shall assume this below. Then according to [Th], $M$ 
belongs to one of the
following three mutually exclusive categories:
\roster
\item"(I)" $M$ is a Seifert-fibred space admitting no 
essential, i.e.,
incompressible and non-$\partial$-parallel torus.
\item"(II)" $M$ is a hyperbolic manifold.
\item"(III)" $M$ contains an essential torus.
\endroster
\par It turns out to be convenient to consider these three 
cases separately. In
case (I) it is well known that one can completely classify 
the finite/cyclic
fillings of $M$. Considering the torus knots for instance, 
one sees that there
exist infinitely many knots whose exteriors are of type 
(I), each of which
admit an infinity of distinct finite (cyclic or noncyclic) 
surgery slopes. Our
contributions deal with the cases (II) and (III). For the 
former we obtain
\thm{Theorem A} Let $K$ be a knot in a closed, connected, 
orientable,
$3$-manifold $W$, such that the interior of 
$M=W-\roman{int}N(K)$ admits a
complete hyperbolic structure of finite volume.
\roster
\item"(1)" There are at most six finite/cyclic surgeries 
on $K$, and
$\Delta(r_1,r_2)\leq 5$ for any two finite/cyclic surgery 
slopes $r_1$ and
$r_2$ of $K$.\endroster\eject
\roster
\item"(2)" If $r_1$ is a finite/cyclic surgery slope of 
$K$ and $r_2$ is a
cyclic surgery slope of $K$, then $\Delta(r_1,r_2)\leq 2$.
\endroster
\ethm
\par The known realizable maximal number of finite/cyclic 
surgeries on a knot
as in Theorem A and their maximal mutual distance is 5 and 
3 [We]. This
example may also be used to show that Theorem A(2) is sharp.
\par To discuss case (III), we introduce the following 
notion. A compact,
connected, orientable, 3-manifold $M$, with boundary a 
torus, is called a {\it
generalized $n$-iterated torus knot exterior} if $M$ can 
be decomposed along
disjoint, essential tori into a union of $n$ cabled spaces 
(in the sense of
[GLi]) and a Seifert-fibred space which has a Seifert 
fibration over the 2-disc
with exactly two exceptional fibres.
\thm{Theorem B} Let $K$ be a knot in a closed, connected, 
orientable,
$3$-manifold $W$, such that $M=W-\roman{int}N(K)$ is 
irreducible and contains
an essential torus.
\par $(1)$ If $M$ is not a generalized $1$-iterated torus 
knot exterior, then
$\Delta(r_1,r_2)\leq 1$ for any two finite/cyclic surgery 
slopes $r_1$ and
$r_2$ of $K$. In particular, there are at most three 
finite/cyclic surgeries on
$K$.
\par $(2)$ If $M$ is not a generalized $1$- or 
$2$-iterated torus knot exterior
and if $K$ admits a cyclic surgery of odd order, then $K$ 
does not admit any
other finite/cyclic surgery.
\ethm
\par Finite/cyclic fillings on a generalized 1- or 
2-iterated torus knot
exterior $M$ can be completely described. This is 
essentially done in [BH,
\S2], where it is shown that if $M$ is not a union of the 
twisted $I$-bundle
over the Klein bottle and a cabled space, then there are 
no more than six
finite/cyclic fillings of maximal mutual distance 5 
(realized on the
exterior of the $(11,2)$-cable over the $(2,3)$-torus knot 
in $S^3)$. In
particular, it is proved that an iterated torus knot in 
$S^3$, admitting a
nontrivial finite surgery, must be a cable over a torus 
knot. A complete list
of all finite surgeries on cabled knots over torus knots 
in $S^3$ is given in
\S2 of that paper.
\par It is shown in [BZ1, Example 10.6] that Theorem B(1) 
is sharp. We also
note that as the finite/cyclic fillings on generalized 1- 
or 2-iterated torus
knot exteriors are readily determined, Theorem B(2) 
completes the
classification of finite/cyclic surgeries on knots in 
manifolds of odd-order,
cyclic fundamental group whose exteriors contain an 
essential torus.
\par Consider now surgery on knots in the 3-sphere $S^3$. 
As is usual, slopes
for a knot in $S^3$ are parameterized by $\bold 
Q\cup\{\tfrac 10\}$ through the
use of the standard meridian-longitude coordinates [R]. In 
$S^3$ only the
trivial knot admits a $\bold Z$-surgery [Ga3].
\thm{Corollary C} Let $K\subset S^3$ be a hyperbolic knot.
\roster\item"(1)" There are at most six finite surgeries 
on $K$, and
$\Delta(r_1,r_2)\leq 5$ for any two finite surgery slopes 
$r_1$ and $r_2$ of
$K$.
\item"(2)" If $r_1$ is a finite surgery slope of $K$ and 
$r_2$ is a cyclic
surgery slope of $K$, then $\Delta(r_1,r_2)\leq 2$. In 
particular, if  $m/n$ is
a finite surgery slope of $K$, then $|n|\leq 2$.
\endroster
\ethm
\par It is shown in [BH] that the $(-2,3,7)$-pretzel knot 
admits at least four
finite surgeries of maximal mutual distance 2. We prove in 
[BZ1, Example
10.1] that this knot has no other finite slopes. This 
example exhibits the
known maximal number of finite surgeries on a hyperbolic 
knot in $S^3$.

\thm{Corollary D} Let $K\subset S^3$ be a satellite knot. 
If $K$ admits a
nontrivial finite surgery, then $K$ is a cabled knot over 
a torus knot.
\ethm
\par It follows from our previous remarks that Corollary D 
classifies finite
surgeries on satellite knots in $S^3$.
\par In their recent work [BH], Bleiler and Hodgson 
obtained, using a
completely different approach, the number 24 and the 
distance 23 for
finite/cyclic surgery on a knot as in Theorem A and the 
number 8 and the
distance 5 for finite/cyclic surgery on a knot as in 
Theorem B.
\par It is a classic result [Mi] that a finite group, 
which is the fundamental
group of a 3-manifold, must belong to one of the following 
types: C-type,
cyclic groups; D-type, dihedral-type groups; T-type, 
tetrahedral-type
groups; O-type, octahedral-type groups; I-type, 
icosahedral-type groups;
and Q-type, quaternionic-type groups. 
It is shown in [BZ1] that more precise
information on finite/cyclic surgeries of a given type may 
be obtained. For
knots in $S^3$ these yield the following result.
\thm{Proposition E} Let $K\subset S^3$ be a hyperbolic knot.
\par{\rm(1)} {\rm(i)} Any $\roman D$-type slope of $K$ 
must be an 
integral slope.
\roster
\item"{\rm(ii)}" There is at most one $\roman D$-type 
finite surgery on $K$.
\item"{\rm(iii)}" If there is a $\roman D$-type surgery on 
$K$, then there is no
even-order cyclic surgery on $K$.
\item"{\rm(iv)}" If there is a $\roman D$-type surgery 
$\alpha$ on $K$, then there is at
most one nontrivial cyclic surgery on $K$; and if there is 
one, $\beta$, say,
then $\alpha$ and $\beta$ are consecutive integers.
\endroster
\par $(2)$ There are at most two $\roman T$-type slopes on 
$K$; and if two, one is
integral, the other has denominator $2$, and their 
distance is $3$.
\par $(3)$ {\rm(i)} Any $\roman O$-type slope of $K$ must 
be an integral slope.
\roster
\item"{\rm(ii)}" There are at most two $\roman O$-type 
slopes on $K$; 
and if two, their
distance is $4$.
\item"{\rm(iii)}" If there is an $\roman O$-type surgery 
on $K$, then there is no
even-order cyclic surgery on $K$.
\item"{\rm(iv)}" If there is an $\roman O$-type surgery 
$\alpha$ on $K$, then there is
at most one nontrivial cyclic surgery on $K$\RM; and if 
there is one, $\beta$,
say, then $\alpha$ and $\beta$ are consecutive integers.
\item"{\rm(v)}" If there is an $\roman O$-type surgery 
slope and a 
$\roman D$-type surgery slope
on $K$, then they are consecutive even integers.
\endroster
\ethm
\par Proposition E may be complemented by various examples 
of hyperbolic knots
in $S^3$ which admit O-type, or I-type, or D-type surgery 
[BH]. In [BZ1]
we produce a few more examples. Notably, [BZ1, Example 
10.2] provides a
hyperbolic knot in $S^3$ which admits a D-type surgery and 
an odd-order
nontrivial cyclic surgery, and [BZ1, Example 10.4] 
provides a hyperbolic knot
in $S^3$ which admits a T-type surgery and a nontrivial 
cyclic surgery.
\par Based on the results obtained in [BZ1] and known 
examples, we raise the
following  problem.
\ex{Finite/cyclic surgery problem} (I) Let $K$ be a knot 
in a connected,
closed,
orientable 3-manifold $W$ such that 
$M=W-\roman{int}(N(K))$ is a hyperbolic
manifold. Show that there are at most five finite/cyclic 
surgeries on $K$ and
that the distance between any two finite/cyclic surgery 
slopes is at most
3.
\par(II) Let $K\subset S^3$ be a hyperbolic knot. Show that:
\roster
\item"(1)" there are at most four finite surgeries on $K$;
\item"(2)" the nontrivial finite surgery slopes on $K$ 
form a set of consecutive
integers;
\item"(3)" the distance between any two finite surgery 
slopes of $K$ is at most
2;
\item"(4)" there is at most one finite surgery on $K$ with 
an even integral
slope.
\endroster
\endex
\par Evidence supporting a positive solution to this 
problem may be found in
[BZ1]. For instance, if the minimal norm [CGLS, Chapter 1] 
amongst all nontrivial
elements of $H_1(\partial M)$ is greater than 24 (16 for 
knots whose exteriors
have no 2-torsion in their homology), then the methods of 
[BZ1] show that for a
knot $K$ as in Theorem A there are at most four 
finite/cyclic surgeries on $K$
of maximal mutual distance no more than 2. These methods 
may also be used to
solve the finite/cyclic surgery problem for various 
families of knots, such as
those with 2-bridges [Ta].
\par Proofs of the results listed above may be found in 
[BZ1]. The techniques
used there are based on the work on M. Culler, C. M. 
Gordon, J. Luecke, and P.
Shalen [CS, CGLS]. In particular, we derive results on the 
norm defined in
Chapter 1 of [CGLS] which, when combined with Chapter 2 of 
that paper, yields
Theorem A. Theorem B is proven by considering the torus 
decomposition of $M$
and then applying Theorem A and results from [Gal, Ga2, 
Sch, CGLS] to determine
what the pieces of this decomposition are, under the added 
assumption that
there are two finite/cyclic slopes such that either (i) 
they are of distance
greater than 1 apart, or (ii) one of them is a cyclic 
slope of odd order.
\par Our study of exceptional fillings is continued in 
[BZ2].
\thm{Theorem F} Let $M$ be a compact, connected, 
orientable, irreducible
$3$-manifold with $\partial M$ a torus. Assume that $M$ is 
not an atoroidal
Seifert-fibred manifold. Fix slopes $r_1$ and $r_2$ on 
$\partial M$ and
suppose that $M(r_1)$ is a reducible manifold.
\par$(1)$ If $M(r_2)$ has a cyclic fundamental group, then
$\Delta(r_1,r_2)\leq 1$.
\par$(2)$ If $M$ is hyperbolic and $r_2$ is a finite 
slope, then
$\Delta(r_1,r_2)\leq 5$ unless $M(r_1)=\bold RP^3\#\bold 
RP^3$ and
$\pi_1(M(r_2))$ is a $D$-type group or a $Q$-type group.
\par$(3)$ If $M$ contains an essential torus and $r_2$ is 
a finite slope, then
$\Delta(r_1,r_2)\leq 1$ unless $M$ is a cable on the 
twisted $I$-bundle over the
Klein bottle or a cable on a hyperbolic manifold for which 
the inequality of
part $(2)$ does not hold.
\ethm
\par We also obtain a new proof of the following result of 
Gordon and Luecke.
\thm{Theorem G {\rm([GLu1])}} Let $M$ be a compact, 
connected, orientable,
irreducible $3$-manifold with $\partial M$ a torus. If 
$M(r_i)$ is a reducible
manifold, for $i=1,2$, then $\Delta(r_1,r_2)\leq 1$.
\ethm
\par Combining the last two results with the cyclic 
surgery theorem, we obtain
\thm{Corollary H} Let $M$ be a compact, connected, 
orientable, irreducible
$3$-manifold with $\partial M$ a torus. Suppose that $M$ 
is not an atoroidal
Seifert-fibred space. If for $i=1,2$, $M(r_i)$ is either a 
reducible manifold
or a manifold with cyclic fundamental group, then 
$\Delta(r_1,r_2)\leq 1$.
Consequently there are at most three cyclic/reducible 
fillings on $M$.
\ethm
\par Suppose now that $r_1$ is a slope on $\partial M$ 
such that $M(r_1)$ has
the fundamental group of a Seifert-fibred manifold $W$ and 
that $r_2$ is a
finite/cyclic slope. Theorems A and B provide upper bounds 
for
$\Delta(r_1,r_2)$ when $W$ admits a Seifert-fibration with 
base orbifold the
$2$-sphere having no more than two exceptional fibres or 
three such fibres if
their indices form a platonic triple. Our next theorem 
deals with most of the
remaining cases.
\thm{Theorem I} Let $M$ be a compact, connected, 
orientable, $3$-manifold with
$\partial M$ a torus. Suppose further that $M$ is neither 
Seifert-fibred nor a
cable on a Seifert-fibred manifold. Let $r_1$ be a slope 
on $\partial M$ such
that $M(r_1)$ has the fundamental group of a 
Seifert-fibred space which admits
no Seifert fibration having base orbifold the $2$-sphere 
with exactly three
exceptional fibres. Then
\par$(1)$ $\Delta(r_1,r_2)\leq 1$ if $M(r_2)$ has a cyclic 
fundamental
group\,\RM;
\par$(2)$ $\Delta(r_1,r_2)\leq 5$ if $M(r_2)$ has a finite 
fundamental group
unless $M(r_1)$ is either $\bold RP^3\#\bold RP^3$ or a 
union of two copies of
the twisted $I$-bundle over the Klein bottle, and 
$\pi_1(M(r_2))$ is a $D$-type
group or a $Q$-type group.
\ethm
\par Applying this result to knots in $S^3$, we obtain
\thm{Corollary J} Let $M$ be the exterior of a knot $K$ in 
$S^3$ and $r$ a
slope on $\partial M$ such that $M(r)$ has the fundamental 
group of a
Seifert-fibred space.
\par$(1)$ If $M(r)$ has the fundamental group of a 
Seifert-fibred space which
is Haken, then $r$ is an integral slope.
\par$(2)$ If $K$ is a satellite knot which is not cabled 
exactly once, then $r$
is an integral slope.
\ethm
\par Finally we return to the question of quantifying 
nonhyperbolic slopes on a
hyperbolic manifold $M$. Call a slope on $\partial M$ {\it 
big Seifert}, or
$bS$ for short, if the associated filling yields a 
Seifert-fibred manifold
whose base orbifold is not a 2-sphere with fewer than four 
cone points. Using
the results above as well as results from [GLi] and [Go2], 
we may prove the
following.
\thm{Corollary K} Let $M$ be a compact, connected, 
orientable, hyperbolic
$3$-manifold with $\partial M$ a torus. Let $r$ and $s$ be 
two slopes that are
contained in the set of all reducible/cyclic/finite/bS 
slopes. If neither
$M(r)$ not $M(s)$ is $\bold RP^3\#\bold RP^3$ or a union 
of two copies of the
twisted $I$-bundle over the Klein bottle, then 
$\Delta(r,s)\leq 5$. Hence there
are at most eight such slopes.
\ethm

\par The reader will find further discussion and results
on exceptional
slopes
in [BZ2].
\heading Acknowledgment\endheading
\par The authors thank Marc Troyanov for several valuable 
conversations
concerning the material in this article.

\enddocument